\documentclass[a4paper,reqno]{amsart}
\usepackage{newcent} 
\usepackage[T1]{fontenc}
\usepackage[centertags]{amsmath}
\usepackage{amsfonts}
\usepackage{amssymb}
\usepackage{amsthm}
\usepackage{newlfont}
\usepackage{fancyhdr,fancyvrb}
\usepackage{graphicx}
\usepackage{pinlabel}
\usepackage{nextpage}
\usepackage{color}
\usepackage{hyperref}
\usepackage{esint}
\usepackage{soul}
\usepackage{Mat-Mor-14-private}

\numberwithin{equation}{section}

\newtheorem{defin}{Definition}[section]
\newtheorem{theorem}[defin]{Theorem}

\theoremstyle{definition} {}

\title[Homogenization problems in the calculus of variations: an overview]{Homogenization problems in the calculus of variations: \\ an overview
 \vspace{.5cm}
 \center{\tiny {dedicated to prof. orlando lopes}}}%
\author{Jos\'e Matias}
\address{CAMGSD, Departamento de Matem\'atica, Instituto Superior T\'ecnico, Av.\@ Rovisco Pais, 1, 1049-001 Lisboa, Portugal}
\email[J.~Matias \myenv]{jose.c.matias@tecnico.ulisboa.pt}
\author{Marco Morandotti}
\address{SISSA -- International School for Advanced Studies, Via Bonomea, 265, 34136 Trieste, Italy} 
\email[M.~Morandotti]{marco.morandotti@sissa.it}

\date{January 29, 2015. Preprint SISSA: 13/2015/MATE}

\makeindex%
\begin{document}

\begin{abstract}
In this note we present a brief overview of variational methods to solve homogenization problems. 
The purpose is to give a first insight on the subject by presenting some fundamental theoretical tools, both classical and modern.
We conclude by mentioning some open problems.
\end{abstract}
\maketitle%
{\small

\bigskip
\bigskip

\tableofcontents%

\section{Introduction}

In this work we intend to give a brief overview on homogenization results derived through variational methods. 
Keeping in mind that the reader is not necessarily familiar with the techniques involved, we will start by introducing the Direct Method of the Calculus of Variations, which combines lower semicontinuity of an energy functional and compactness properties of minimizing sequences to grant the existence of a minimum for the energy.
Minimizers of an energy functional describing the state of a physical system are sought for since they describe the equilibrium configurations of the system.
Moreover, it will be necessary to study the limit behaviour of a family or of a sequence of energy functionals of integral type,  and results that grant that the limit functional admits an integral representation as well will be most useful.

The mathematical theory of homogenization was motivated from problems in Physics and Continuum Mechanics, in order to describe the behavior of composite materials and reticulated structures. 
The former are characterized by different constituents finely mixed in a structured way that bestows enhanced properties on the composite material. 
The latter are characterized by periodically distributed holes or inclusions. 
In both cases, the material can be modeled as a grid-like structure in which every cell is the copy of an elementary one which contains both materials. 
Layered materials are another type of periodic composites where the repeating structure is not a cell, but a whole layer.
The way the individual properties of the components determine the macroscopic behavior of the composite material is far from trivial and it is the object of the mathematical theory of homogenization.
This theory provides the tools for deducing the relevant properties of the compound by averaging out (in a sense that needs to be precisely stated) the properties and the proportions of its components. 
At the modeling stage, the fine structure of the composite is introduced by means of a smallness parameter $\eps>0$, which appears in the energy functional. Homogenizing the functional means to study the limit behavior as $\eps$ tends to $0$, and to obtain a new functional that does not explicitly account for the microstructure. The limiting procedure is a very delicate one and it is (usually) performed by means of $\Gamma$-convergence \cite{DGF,DG,DM,B2}.
The so-called $\Gamma$-limit might have a very different structure from that of the $\eps$-functionals (see, e.g., \cite{MO}), hence the necessity of the integral representation theorems mentioned above.
The limit energy density is usually found by solving a minimum problem in the elementary cell and computing averages over larger and larger cells. 

The overall plan of this work is as follows: we close the introduction by establishing the notation we are going to use throughout the paper; in Section \ref{sec:dirmet} we give a brief overview of the Direct Method of the Calculus of Variations; in Section \ref{sec:hom} we explain briefly what a homogenization problem is, give a brief introduction to $\Gamma$-convergence and list some integral representation results which are often used in this context; in Section \ref{sec:res} we state some of the most important homogenization results in the context of Calculus of Variations, distinguishing between the cases of superlinear growth and linear growth of the energy density.
We will close this note by referring, in Section \ref{sec:persp}, to some other techniques used to tackle homogenization problems and we mention some perspective work.

\subsection{Notation}
Throughout the text  we will use the following notations:
\begin{itemize}
\item[-] $\Omega \subset \R{N}$ denotes a bounded open set;
\item[-]  $\cL^{N}$ is the $N$-dimensional Lebesgue measure. If $A$ is a Lebesgue measurable set,
$\cL^{N}(A)$ denotes its measure. Sometimes also the notation $|A|$ will be used;
\item[-] $\cS^{d-1}$ denotes the unit sphere in $\R{d}$;
\item[-] $Q:=[-\frac12,\frac12]^N$ denotes the unit cube in $\R{N}$.
\item[-] $\cO(\Omega)$ denotes the family of the open subset contained in $\Omega$.
\item[-] $C^k(\Omega;\R{d})$ denotes the space of $\R{d}$-valued functions which have continuous derivatives up to order $k$.
\item[-] $\Lp1{Q-\per}(\R{N};\R{d})$ denotes the space of $\R{d}$-valued $L^1$ functions defined on $Q$ and extended by $Q$-periodicity to the whole of $\R{N}$; 
\item[-] $\cM(\Omega;\R{d})$ denotes the space of $\R{d}$-valued Radon measures defined on $\Omega$;
\item[-] $\cD'(\Omega;\R{d})$ denotes the space of $\R{d}$-valued distributions defined on $\Omega$;
\item[-] $W^{1,\infty}(\Omega;\R{d})$ denotes the space of $\R{d}$ valued functions which, as well as their weak first order partial derivatives, are essentially bounded;
\item[-] $BV(\Omega;\R{d})$ denotes the space of $\R{d}$-valued functions with bounded variation defined on $\Omega$;
\item[-] $S_u$ denotes the approximate discontinuity set of $u\in BV(\Omega;\R{d})$;
\item[-] $[u](x)$ denotes the jump of $u$ across $x\in S_u$.
\item[-] $C$ represents a generic positive constant, which may vary from expression to expression;
\end{itemize}

\section{The Direct Method of the Calculus of Variations}\label{sec:dirmet}
The origins of the direct method in the Calculus of Variations goes back to Hilbert, for treating the Dirichlet integral, and to Lebesgue and Tonelli. 
Let $X$ be a normed space and let $I:X\to[-\infty, +\infty]$ be a functional which is not identically $\infty$.
The direct method developed by Tonelli states which conditions $X$ and $I$ have to satisfy to have the existence of a minimum point for $I$, that is the existence of $\bar u\in X$ such that $I(\bar u)=\inf_{u\in X} I(u)$.
Following the enlightening description given in \cite{FonLeo}, the direct method consists in the recipe:
\begin{itemize}
\item[(DM1)] Consider a \emph{minimizing sequence} $\{u_n\}\subset X$, that is, a sequence that converges to the infimum of $I$:
$$\lim_{n\to\infty} I(u_n)=\inf_{u\in X} I(u);$$
\item[(DM2)] Prove that $\{u_n\}$ admits a subsequence $\{u_{n_k}\}$ that converges with respect to some topology $\tau$ to some point $\bar u\in X$;
\item[(DM3)] Establish the sequential lower semicontinuity of $I$ with respect to $\tau$;
\item[(DM4)] Conclude that $\bar u$ is a minimum of $I$; indeed:
$$\inf_{u\in X} I(u)=\lim_{n\to\infty} I(u_n)=\lim_{k\to\infty} I(u_{n_k})\geq I(\bar u)\geq\inf_{u\in X} I(u).$$
\end{itemize}
The step (DM2) usually can be done by coercivity and growth conditions on $I$, while (DM3) is the challenging part in the program.

Another important observation is that (DM3) indeed fails in most applications, hence there is the need to overcome this impasse.
This problem is central in the calculus of variations, and is usually achieved by introducing the \emph{relaxed functional} $\tilde I$of $I$.
Such a functional is studied instead of $I$ and is expected to return the same minimum.
For this to happen, $\tilde I$ must enjoy the following properties:
\begin{itemize}
\item[(R1)] $\tilde I$ is sequentially lower semicontinuous with respect to $\tau$;
\item[(R2)] $\tilde I$ inherits coercivity from $I$;
\item[(R3)] $\min_{u\in X} \tilde I(u)=\inf_{u\in X} I(u)$.
\end{itemize}

One of the most important classes of functionals in the calculus of variations is given by integral functionals, that is, functionals of the form
\begin{equation}\label{intfun}
I(u)=\int_{\Omega} f(x,u(x),\nabla u(x),\ldots)\,\de x,
\end{equation}
where $u:\Omega\to\R{d}$, $d\geq1$, and $f:\Omega\times\R{d}\times\R{d\times N}\times\cdots\to [-\infty,\infty]$ is usually referred to as \emph{lagrangian} or \emph{energy density}.
Which derivatives $f$ depends on varies from problem to problem.
The simplest cases, for which a sound theory is available, are dependences up to the first gradient of $u$. 
The following result, stated in the context of Sobolev spaces, (see \cite[Theorem 3.3]{DacIntro}) contains all the general ideas of the direct method.
\begin{theorem}
Let $\Omega\subset\R{N}$ be a bounded open set with Lipschitz boundary. Let $f\in C^0(\overline\Omega\times\R{}\times\R{N})$, $f=f(x,u,\xi)$, satisfy
\begin{itemize}
\item[(H1)] $\xi\to f(x,u,\xi)$ is convex for every $(x,u)\in\overline\Omega\times\R{}$;
\item[(H2)] there exist $p>q\geq1$ and $\alpha_1>0$, $\alpha_2,\alpha_3\in\R{}$ such that
$$f(x,u,\xi)\geq \alpha_1|\xi|^p+\alpha_2|u|^q+\alpha_3,\qquad\text{for all $(x,u,\xi)\in\overline\Omega\times\R{}\times\R{N}$}.$$
\end{itemize}
Let 
\begin{equation}\label{minprob}\tag{P}
\inf\left\{ I(u)=\int_\Omega f(x,u(x),\nabla u(x))\,\de x:\;u\in u_0+W^{1,p}_0(\Omega)\right\}=m
\end{equation}
where $u_0\in W^{1,p}(\Omega)$ with $I(u_0)<\infty$. 
Then there exists $\bar u\in u_0+W^{1,p}_0(\Omega)$ a minimizer of \eqref{minprob}.

Furthermore, if $(u,\xi)\to f(x,u,\xi)$ is strictly convex for every $x\in\overline\Omega$, then the minimizer is unique.
\end{theorem}

Hypothesis (H2) ensures (DM2), (H1) ensures (DM3), therefore the direct method gives the existence of a minimizer, while strict convexity is usually the requirement needed to prove that the minimizer is unique.

In the case where $u:\Omega \subset \R{N} \to\R{d}$ with $N, d > 1$, convexity of the integrand $f$ in the last variable is still sufficient but it is no longer a necessary condition for sequential lower semicontinuity.
The notion of \emph{quasiconvexity} introduced by Morrey \cite{MO} is the generalization of convexity to this case (usually referred to, in the literature, as the vectorial case), and it is a necessary and sufficient condition to the weak sequential lower semicontinuity of $I$ on Sobolev spaces.

A function $f$ is said to be quasiconvex at $\xi\in \R{N{\times}d}$  and 
\begin{equation}\label{qc}
f(\xi)\leq \int_Q f(\xi+ \nabla \Phi(x))\,\de x, \quad\text{for all $\Phi\in C^\infty_c(Q; \R{d})$,}
\end{equation}

In order for $f$ to be quasiconvex at $\xi\in\R{N{\times}d}$, its value at $\xi$ has to be minimum with respect to averaged variations by means of gradients of smooth functions.
Condition \eqref{qc} can be written in the following way
\begin{equation}\label{qc1}
f(\xi)\leq \int_Q f(\xi+v(x))\,\de x, \quad\text{for all $v\in C^\infty_c(Q; \R{d})$, such that $\curl v=0$},
\end{equation}
since a vanishing $\curl$ is the condition for a function to be a gradient, $v=\nabla\Phi$.
Moreover, if $f$ is quasiconvex, then one can show \cite{FM} that the following (variational) characterization holds
\be \label{qc2}
f(\xi) = \inf_{\Phi \in W^{1,\infty}_{\text {per}} (Q; \R{d})}\int_Q f(\xi + \nabla \Phi(x))\; \de x,
\ee
where $W^{1,\infty}_{\text{per}}(Q; \R{d})$ stands for the class of periodic functions in $W^{1, \infty} (Q; \R{d}).$
It must be noted (see \cite{T}) that many vector-valued problems coming from applications to continuum mechanics and electromagnetism are modeled by integral functionals whose integrand satisfies \eqref{qc1}, but with operators $\cA$ different than $\curl$: one seeks for the property, called \emph{$\cA$-quasiconvexity},
\begin{equation}\label{Aqc}
f(\xi)\leq \int_Q f(\xi+v(x))\,\de x, \quad\text{for all $v\in C^{\infty}_{\text{per}}(Q; \R{d})$, such that $\cA v=0$ and $\int_Q v(x)\,\de x=0$}.
\end{equation}

The notion of $\cA$-quasiconvexity was first introduced by Dacorogna \cite{Dac} following the works of Murat and Tartar 
in compensated compacteness (see \cite{Mur_84} and \cite{T}) and was further
devoloped by Fonseca and M\"{u}ller \cite{FM}  (see also \cite{BFL}).


In \eqref{Aqc}, $\cA: \cD'(\Omega;\R{d}) \to  \cD'(\Omega;\R{M})$ is a linear first order partial differential operator with constant coefficients 
of the form
\be\label{operator} 
\cA = \sum_{i=1}^N A^{(i)} \frac{\partial}{\partial x_i},
\ee
where $A^{(i)} \in \M{M{\times} d}, M\in\N{}$, are called the \emph{symbols} of the operator.
Also, $\cA$ is assumed to satisfy the 
 so-called Murat's condition of {\it constant rank} (introduced in \cite{Mur_84}; see also \cite{FM}), i.e., there exists $c\in \N{}$ such that 
\be\label{murat}
\mathrm{rank}\, \left(\sum_{i=1}^N A^{(i)} \xi_i\right)=c \quad\, \text{for all}\,\,\, \xi=(\xi_1,...,\xi_d)\in \cS^{d-1}.
\ee
This condition is of technical nature and is necessary in order to define a projection operator in the kernel of $\cA$ (see \cite{FM}) and to obtain technical lemmas. We mention it here because the results cited in the rest of the paper are obtained using \eqref{murat}.

We now present examples of such operators (for more examples see \cite{FM}):
\begin{enumerate}
\item  ($\cA= \text{div}$) For $\mu \in {\cM}(\Omega; \R{N})$ we define
$$\cA \mu = \sum_{i=1}^N \frac {\partial \mu^i} {\partial x_i}.$$
\item ($\cA=\text{curl}$) For  $\mu \in {\cal M}(\Omega; \R{N{\times}m})$,  we define
$$\cA \mu = {\left( \frac {\partial \mu_k^j} {\partial x_i} -  \frac {\partial \mu_i^j} {\partial x_k} \right) }_{j=1,..,m;\, i,k=1,..,N}.$$
\item (Maxwell's Equations)  For  $\mu= (m,h) \in {\cM}(\R{3}; \R{3 \times 2})$ we define
$$ \cA \mu = \left( \text{div} (m+h), \text{curl} \,h \right).$$
\end{enumerate}
All the examples above are structured in a way that it will be easy later on to impose $\cA\mu=0$ as a constraint to a minimum problem.

We close this section by giving the formal definition of $\cA$-quasiconvexity and by introducing two objects that are naturally attached to the operator $\cA$.
\begin{defin}[$\cA$-quasiconvex function]\label{aquasi} 
A locally bounded Borel function $f:\R{d}\to \R{}$ is said to be  $\cA$-quasiconvex if
\[
f(v)\leq \int_Q f(v+w(x))\, \de x
\]
 for all $v\in \R{d}$ and for all  $w\in C_{Q-\per}^\infty(\R{N};\R{d})$
such that $\cA w=0$ in $\R{M}$ with $\int_Q w(x)\, \de x=0$.
\end{defin}
We remind the reader that in the case $\cA = \text{curl}$, $\cA$-quasiconvexity reduces to quasiconvexity.

We conclude this section by introducing two objects that are naturally attached to the operator $\cA$. 
\begin{defin}
\label{355}
(i) The \emph{characteristic cone} of $\cA$ is defined by
\be\label{356}
\cC:=\left\{v\in\R{d}:\exists w\in\R{N}\setminus\{0\}, \left(\sum_{i=1}^N A^{(i)}w_i\right)v=0\right\}.
\ee
(ii) Given $v\in\R{d}$, we define the linear subspace of $\R{N}$
$$
\cV_v:=\left\{w\in\R{N}:\left(\sum_{i=1}^N A^{(i)}w_i\right)v=0\right\}.
$$
If $v\notin\cC$, then $\cV_v=\{0\}$, otherwise $\cV_v$ is a non trivial subspace of $\R{N}$.
\end{defin}

\section{Homogenization}\label{sec:hom}
The mathematical theory of homogenization is a powerful tool to deduce averaged properties of composite materials with a periodic structure.
In each periodicity cell, two or more components with different physical or chemical properties coexist, and their interplay can be tuned to tailor the material to enhance some specific characteristics.
The general idea is that of finding a way to disregard the periodic structure and to look at the material at a larger (length-)scale.
The desired macroscopic properties, that will be the result of the microscopic structure, are usually captured via some limiting process.
One should imagine to have a composed material and to look at it from far away in such a way not to be able to distinguish the fine details anymore, but to grasp the effect of their combination.

A homogenization problem can be mathematically described in essentially two ways: one considering a PDE, the other one considering a functional. 
Both the PDE and the functionals are indexed by a smallness parameter $\eps>0$ and the scope is to study the limit -- in some appropriate sense -- either of the PDE or of the functional, as $\eps\to0$.

The approach via functional minimization is the one proper of the calculus of variations and the right notion of convergence is that of $\Gamma$-convergence, which we briefly describe in Subsection \ref{gammac}.

When the functionals involved are of integral type, e.g., as in \eqref{intfun}, one of the main challenges is to prove that also the limiting functional given by $\Gamma$-convergence, the $\Gamma$-limit, is of integral type too. 
This is a central question in homogenization theory, since $\Gamma$-limits fo functionals can have a sensibly different form compared to the $\eps$-functionals (see, for example, the case of the Cahn-Hilliard functional for fluid mixtures: at the level $\eps$ the functionals involve a Dirichlet energy plus a potential term, whereas their $\Gamma$-limit is the perimeter functional \cite{MO}).
The issue of integral representation of functionals will be addressed in Subsection \ref{intrep}

\subsection{$\Gamma$-convergence}\label{gammac}
The scope of $\Gamma$-convergence is the description of the asymptotic behavior of families of minimum problems of the type
\be \label{hom1} \text{min} \{ F_\eps (u): u \in X_\eps \}
\ee
by finding a limiting problem 
\be \label{hom2} \text{min} \{\mathcal{F}(u): u \in X \}
\ee
whose solution captures the relevant behavior of the minimizers of $F_\eps$.
In this way, the main properties of the solutions of \eqref{hom1} can be approximately described by those of the solutions of \eqref{hom2}. 

Before stating the main result on $\Gamma$-convergence we briefly describe the main ideas behind this powerful technique introduced by De Giorgi (\cite{DG}, \cite{DGF}). 
This description follows closely the one given in the introduction of \cite{B2}.

A preliminary step in order to apply $\Gamma$-convergence is to ensure compactness, that is, ensure that sequences of minimizers $u_\eps$ of \eqref{hom1}, if they exist, converge in some sense to some $u \in X$. 
The candidate space $X$ for the limit problem \eqref{hom2} is then determined by this compactness argument. 
Then, the functional $\mathcal{F}$ in the limit problem \eqref{hom2} is determined via a process of optimization between lower and upper bounds. 
A lower bound for $\mathcal{F}$ is an energy $G$ that satisfies
\be \label{lbound} G(u) \leq \liminf _{ \eps \to 0^+} F_\eps (u_\eps), \; \text{whenever }\; u_\eps \to u.
\ee
This requirement, implies some structure properties for the candidate $G$, in particular lower semicontinuity. 
Condition \eqref{lbound} implies that
\be \label{inf}
\inf\{ G(u): u \in X\} \leq \lim_{\eps \to 0^+} \text{min} \{ F_\eps (u): u \in X_\eps\},
\ee
given this limit exists. 
The sharpest lower bound is then derived by optimizing the role of $G$ and this is done through a minimization argument. 
The derivation of the optimal $G$ suggests the form of the minimizing sequences. 
We need now to construct, for each $u \in X$, a particular sequence $\overline{u}_\eps \to u$ and define $H(u) = \lim_{\eps \to 0^+} F_\eps(\overline{u}_\eps).$ This functional $H$ is an upper bound for the limit energy, and for such $H$ we then have that:
\be \label{upper}
\lim_{\eps \to 0^+} \text{min} \{ F_\eps (u): u \in X_\eps\} \leq \inf \{ H(u): u \in X\}.
\ee
If one can prove that these two bounds coincide, then this is the $\Gamma$-convergence of $F_\eps$ in \eqref{hom1} to $\cF$ in \eqref{hom2}.
Having defined the upper and lower bounds on the energies for all functions an not only for minimizers, $\Gamma$-convergence implies the following useful properties: (a) convergence of minimum problems \eqref{hom1} to the minimum problem \eqref{hom2}, and of minimizers, (b) stability of $\Gamma$-convergence under continuous perturbations, and (3) lower semicontinuity of the $\Gamma$-limit $\cF$.


We state now the precise definition and main properties of $\Gamma$-convergence and we refer to 
\cite{DM93} for a comprehensive treatment and bibliography on this subject. 
Let $X$ denote a metric space.
\begin{defin}{\rm ($\Gamma$-convergence of a sequence of functionals)} Let
 $F_n ,F:X \rightarrow {\R{}}\cup \{{+\infty}\}$. 
 The functional
${ F} $ is said to be the
$\Gamma$-$\liminf$ (resp.\  $\Gamma$-$\limsup${\rm )} of
$\{{ F}_n\}_{n} $ with respect to the metric of $X$
 if for every $u \in X$

$${ F}(u)= \inf_{\{u_n\}} \left\{
\liminf_{n \to +\infty} { F}_n (u_n): ~ u_n \in
X ,\,  u_n \to u \text{ in } X \right\}\,\, ({\text resp.\ }
\limsup_{n \to +\infty}).
$$
In this case we write
$$ F= \Gamma \text{-}
\liminf_{n \to +\infty}{F}_{n} \enspace
\enspace \left( {\text resp.\ }\,\,\, F= \Gamma
\text{-}\limsup_ {n \to +\infty}{F}_{n}\right).
$$
Moreover,   $F$  is  said to be
the $\Gamma$-$\lim$  of $\{{F}_{n}\}_{n} $  if
$$ F =\Gamma \text{-} \liminf_{n\to +\infty} {F}_{n}=\Gamma
\text{-}\limsup_{n\to +\infty} {F}_{n},$$
and in this case we write
$$ {F}= \Gamma
\text{-}\lim_{n \to +\infty}{F}_{n}.
$$
\end{defin}

For every $\eps>0$ let $F_\eps$ be a functional defined in $X$ with
values in ${\R{}}\cup \{{+\infty}\}$. 

\begin{defin}\label{defgl}{\rm ($\Gamma$-convergence of a family of
functionals)}
A functional  ${F} :X \to {\R{}}\cup \{{+\infty}\}$ is said to be
the $\Gamma \text{-}\liminf$ (resp.\ $\Gamma$-$\limsup$ or
$\Gamma$-$\lim$) of $\{{F}_{\eps}\}_{\eps}$ with respect to
the metric of $X$, as $\eps\to 0^+$, if for every sequence
$ \eps_n \rightarrow 0^+$,
$$ {F}= \Gamma \text{-}
\liminf_{n \to +\infty}{F}_{\eps_n} \quad \left(
{\text resp.\ }\,\,\, {F}= \Gamma
\text{-}\limsup_{n\to +\infty }{
F}_{\eps_n} \,\,\, \hbox{or}  \,\,\, {F}= \Gamma
\text{-}\lim_{n\to +\infty }{
F}_{\eps_n} \right),
$$
\noindent and we write
$$ {F}= \Gamma \text{-}
\liminf_{\eps \to 0^+}{F}_{\eps} \quad \left( {\text
resp.\ }\,\,\,  {F}= \Gamma
\text{-}\limsup_{\eps \to 0^+ }{
F}_{\eps}\,\,\, \hbox{or} \,\,\, {F}= \Gamma
\text{-}\lim_{\eps \to 0^+ }{
F}_{\eps}\right).
$$
\end{defin}

One of the most important properties of
$\Gamma$-convergence is that under appropriate compactness assumptions
it implies the convergence of minimizers of a family of
functionals to the minimum of the limiting
functional, as a consequence of the following result 
(see Corollary~7.20 in \cite{DM93}).
\begin{theorem}{\rm (Fundamental Theorem of $\Gamma$-convergence)}\label{min}
Let  $\{{F}_{\eps}\}_{\eps} $ be a family of 
functionals defined in $X$ and let  
$${F} = \Gamma\text{-} \lim_{\eps
\to 0^+}{F}_{\eps}.$$ 
If $u_\eps$ is a minimizer of  $F_{\eps} $ in $X$ and $u_\eps \to u$ in $X$
then $u$ is a minimizer of $F$ in $X$ and
$$F(u)= \lim_{\eps \to 0^+}F_{\eps}(u_\eps).$$
\end{theorem}

In the context of homogenization, one is interested in encoding in the solutions of \eqref{hom2} the highly oscillatory nature of the solutions of \eqref{hom1}. 
It is therefore natural to minimize among periodic functions to derive the optimal functional $G$ to satisfy \eqref{inf}, as it will be clear in the rest of this note.

One of the great advantages of $\Gamma$-convergence is that it allows to study a (difficult) problem of the type \eqref{hom2} by approximating it by (possibly easier) problems as in \eqref{hom1}. 
In homogenization problems in particular, it suggests that minimizers of \eqref{hom1} oscillate close to their limit following an energetically-optimal locally-periodic pattern.

\subsection{Integral representation of functionals}\label{intrep}
We already noticed that in most applications the system at study is described by an energy functional of the type \eqref{intfun}. 
In this section we state two important results that are often used to derive the integral representation of the limit functional.

{\bf A global method for relaxation \cite{BFM}}:
This work deals with the identification of the integral representation of a broad class of functionals defined in $BV(\Omega;\R{d})\times \cO(\Omega)$ (for a comprehensive introduction to $BV$ spaces see \cite{AmbrosioFuscoPallara00}).
More precisely, the authors derive an integral representation for functionals within the following abstract framework:
$$ \mathcal{F}: BV (\Omega; \R{d})\times \cO(\Omega) \to [0, +\infty[,$$
satisfying
\begin{itemize}
\item $\mathcal{F}(u,.)$ is the restriction to $\cO(\Omega)$ of a Radon measure;
\item $\mathcal{F}(.; A)$ is $L^1$ lower semicontinuous;
\item $\exists\, C > 0$ such that $0 \leq \mathcal{F}(u;A) \leq C\left[ \int_A ( 1 + |\nabla u|(x))\;\de x + |D^su|(A)\right].$
\end{itemize}

The main idea is to show that, for sets of small size,  $\mathcal{F}(u,A)$ behaves like the set function $m(u; A)$ defined in $\mathcal{A_\infty}(\Omega)$ (the family of all Lipschitz subdomains of $\Omega$) by:
$$ m(u; A) := \inf \{ \mathcal{F}(v; A), \; \; v\lfloor \partial A = u\lfloor \partial A, \; u \in BV(\Omega;\R{d})\}.$$

This general framework allows for different applications, namely for phase transition, fracture mechanics, plasticity and image segmentation problems. In particular, it can be applied to the $\Gamma(L^1)$-limit of a family of functionals satisfying appropriate growth conditions. In this context, an application to homogenization problems in $SBV$ involving bulk and interfacial free energies is given, extending to the case $p=1$  results previously obtained in \cite{BDV} for $p>1$.

For $\eps > 0$ and $A \in \cO(\Omega)$, let $\cF_\eps (.; A)$ be defined in $BV(\Omega; \R{d})$ by:
$$ \cF_\eps (u; A) := \left\{ 
\begin{array}{ll} 
\displaystyle \int_A f\left(\frac{x}{\eps}, \nabla u(x)\right)\; \de x + \int_{S_u \cap A} g\left(\frac{x}{\eps} , [u](x), \nu_u(x)\right)\; \de\mathcal{H}^{N-1}, & u \in BV(\Omega; \R{d}),\\
+ \infty, & \text{otherwise},\\
\end{array}
\right.$$
where $\nu_u(x)$ is the normal to $S_u$ at $x$, and $f$ and $g$ satisfy the following conditions:

\begin{itemize}
\item[(A1)] $f: \R{N}\times \R{d\times N} \to [0, +\infty[$ is a Borel function, $Q$-periodic in the first argument, such that
$$ \frac{1}{C}|\zeta| \leq f(x, \zeta) \leq C( 1 + |\zeta|),$$
\noindent for all $\zeta \in \R{d\times N}, $ for all $ x \in \R{N},$ and for some constant $C > 0$;
\item[(A2)] there exist $m \in ]0,1[,\;  L >0$ such that
$$\left | f^\infty (x, \zeta) - \frac{f(x, t\zeta)}{t}\right| \leq \frac{C}{t^m},$$
for all $\zeta \in \R{d\times N}, with ||\zeta|| = 1$, for $t > L$, and for all $x \in \R{N}$; here $f^\infty$ denotes the recession function of $f$ defined by
$$f^\infty(x,\zeta):=\limsup_{t\to+\infty} \frac{f(x,t\zeta)}{t}.$$
\item[(A3)] $g: \R{N} \times \R{d}\times \mathcal{S}^{N-1} \to [0, +\infty[$ is a Borel function, $Q$-periodic in the first argument, satisfying
$$ \frac{1}{C}|\lambda| \leq g(x, \lambda, \nu) \leq C|\lambda|,$$
\noindent for all $x \in \R{N}, \lambda \in \R{d}$, and $\nu \in \mathcal{S}^{N-1}$;
\item[(A4)] There exist $ \alpha \in ]0,1[, \; l > 0$ such that
$$ \left| \overline{g} (x, \lambda, \nu) -\frac{g(x, t\lambda, \nu)}{t}\right| \leq Ct^\alpha,$$
\noindent where $\overline{g}$ is defined by:
$$ \overline{g}(x, \lambda, \nu):= \limsup_{t \to 0} \frac{g(x, t\lambda, \nu)}{t}.$$
\end{itemize}

The integral representation of $\mathcal{F}$, the $\Gamma(L^1)$-limit is derived, and it is shown that:
\begin{theorem}\label{bfm}
Under hypothesis (A1)--(A4), $\mathcal{F} = \mathcal{F}_{\hom},$ which is given by:
$$ \mathcal{F}_{\hom} ( u; A) := \int_A f_{\hom}(\nabla u(x))\; \de x + \int_{S_u \cap A} g_{\hom}([u], \nu_u)\;\de\mathcal{H}^{N-1} + \int_A g^\infty_{\hom}\left(\frac{\de C(u)}{\de |C(u)|}\right)\; \de |C(u)|,$$
where 
$$f_{\hom}(\zeta) := \lim_{T\to +\infty}\frac1{T^N} \inf_{\substack{u \in SBV(TQ; \R{d})\\ u=\zeta \text{ on }\partial (TQ)}} \left\{ \int_{TQ} f(x, \nabla u(x))\; \de x + \int_{S_u \cap TQ} \overline{g}(x, [u], \nu_u)\; \de \mathcal{H}^{N-1}\right\},$$
$$ g_{\hom}(\lambda, \nu) := \lim_{T\to + \infty} \frac{1}{T^{N-1}} \inf_{\substack{u \in SBV(TQ_\nu; \R{d})\\ u= u_{\lambda,\nu} \text{ on } \partial (TQ_\nu)}} \left\{ \int_{TQ_\nu} f^\infty(x, \nabla u)\; \de x + \int_{S_u \cap TQ_\nu} g(x, [u], \nu_u)\; \de \mathcal H^{N-1} \right \},$$
and where 
$$u_{\lambda, \nu}(y) := \left\{ \begin{array}{ll} \lambda & \text{if} \; y\cdot \nu > 0\\
&\\
0 & \text{otherwise}.
\end{array}
\right.
$$
\end{theorem}
{\bf A global method for relaxation in $W^{1,p}$ and in $SBV_p$ \cite{BFLM}}:
This work extends the result in \cite{BFM} to functionals defined in $SBV_p(\Omega; \R{d}) := \{ u \in SBV (\Omega; \R{d}): \nabla u \in L^p(\Omega; \R{d\times N}, \; \mathcal{H}^{N-1}(S_u \cap \Omega) < + \infty \} $ and with superlinear growth.  
This result was proven previously in \cite{BCP}, under additional regularity hypothesis.

Considered a functional
$$\mathcal{F}: SBV_p(\Omega; \R{d})\times \cO(\Omega) \to [0,+\infty[,$$
satisfying, for every $(u, A) \in SBV_p(\Omega; \R{d}) \times \cO(\Omega)$,
\begin{itemize}
\item [(B1)] $\mathcal{F}(u; .)$ is the restriction to $\cO(\Omega)$ of a Radon measure;
\item[(B2)] $\mathcal{F}(u; A) = \mathcal{F}(v; A)$ whenever $u = v$ $\mathcal{L}^N$-a.e. on $A \in \cO(\Omega)$;
\item [(B3)] $\mathcal{F}(. ; A) $ is $L^1(\Omega; \R{d})$ lower semicontinuous;
\item[(B4)] there exists $C > 0$ such that
$$ \frac{1}{C} \left(\int_A |\nabla u|^p\, \de x +
\int_{S_u \cap A} \!\! (1 + [u])\, \de \mathcal{H}^{N-1} \right) 
\leq \mathcal{F}(u; A)
\leq C\left( \int_A(1 + |\nabla u|^p )\, \de x + \int_{S_u \cap A} \!\!( 1 + [u])\, \de \mathcal{H}^{N-1} \right).$$
\end{itemize}
As in \cite{BFM}, the densities of the integral representation of $\mathcal{F}$ are obtained trough the set function
$$ m(u; A) := \inf \{ \mathcal{F}(v; A), \; \;  v \in SBV_p (A; \R{d}), v = u \; \text {on a neighborhood of } \partial A \}$$

The main result reads as follows
\begin{theorem}\label{mainBFLM}
Under hypotheses (B1)--(B4), for every $(u,A) \in SBV_p(\Omega; \R{d})\times \cO(\Omega),$
$$ \mathcal{F}(u; A) = \int_A f(x, u, \nabla u)\; \de x + \int_{S_u \cap A} g(x, u^+, u^-, \nu_u)\; \de\mathcal{H}^{N-1},$$
\noindent where
$$ f(x_0, u_0, \zeta) := \limsup_{\eps \to 0^+} \frac{(m(u_0 + \zeta (\cdot - x_0); Q_\nu(x_0, \eps))}{\eps^N},$$
$$ g(x_0, a, b, \nu) := \limsup_{\eps \to 0^+} \frac{m(u_{x_0, a, b, \nu}; Q_\nu(x_0, \eps))}{\eps_n^{N-1}},$$
for all $x_0 \in \Omega$, $u_0, a, b \in \R{N}$, $\zeta \in \R{d\times N}$, $\nu \in \mathcal{S}^{d-1}$, and where
$$ u_{x_0, a, b, \nu}(x) := \left\{ 
\begin{array}{ll} a, & \text{if } (x - x_0)\cdot \nu > 0,\\
b, & \text{if } (x - x_0)\cdot \nu \leq  0.\\
\end{array}
\right.
$$
Here, $Q_\nu$ denotes the unit cube with two faces perpendicular to $\nu\in\cS^1$.
\end{theorem}
\section{Homogenization results}\label{sec:res}
In this section we give a brief overview of the most relevant homogenization results in the context of Calculus of Variations. 
We distinguish between the cases where the integrand $f$ satisfies a  superlinear growth condition (with exponent $p >1$) and a linear growth condition ($p=1$), since the technical difficulties encountered and the methods are different. 
This is due mainly to the fact that $L^1$ is not a reflexive space and sequences of $L^1$ functions can converge to measures.
This is another way of phrasing that linear growth conditions allow for concentrations.

\subsection{Case $p>1$}
In \cite{FM} it was shown that $\mathcal{A}$-quasiconvexity is a necessary and sufficient condition for sequential lower semicontinuity of functionals of the form
\be \label{acqfunctional} (u, v) \to \int_\Omega f(x, u(x), v(x))\; \de x,\ee
\noindent where $f: \Omega \times \R{m}\times \R{d} \to [0,+\infty[$ is a  normal integrand (that is, $f$ is Borel measurable and $v\mapsto f(x, u, v) $ is lower semicontinuous for a.e.~$(x, u)\in\Omega\times\R{m}$),  $u_n \to u$ in measure, $v_n \wto v $ in $L^p(\Omega; \R{d})$ ($v_n \wsto v$  if $p=+\infty$) and $\mathcal{A}v_n \wto 0$ in $W^{-1,p}$ ($\mathcal{A}v_n = 0$ if $p=+\infty$ ).

Following this result, \cite{BFL}  deals with the integral representation of relaxed energies and of $\Gamma$-limits of functionals of the type \eqref{acqfunctional} under the following hypotheses:
 $f: \Omega \times \R{m} \times \R{d} \to [0,+\infty[$ is a Carath\'{e}odory integrand satisfying
$$ 0 \leq f(x, u, v) \leq a(x,u)(1 + |v|^p),$$
for a.e.\@ $x \in \Omega$ and all $(u,v) \in \R{m}\times \R{d}$, where $ 1 \leq p < +\infty$, $a \in L^{\infty}_{\loc}(\Omega \times \R{}; [0, +\infty[)$, $u_n \to u$ in measure, $v_n \wto v $ in $L^p(\Omega; \R{d})$ and $\mathcal{A}v_n \to 0$ in $W^{-1, p}(\Omega; \R{M})$.

As an application, a  homogenization result  is derived, in the context of $\cA$-quasiconvexity for integrands with growth of order $ p >1$, with one microscopic scale. 
Precisely, the authors present a homogenization result for periodic integrands in the context of $\cA$-quasiconvexity. 
For $\eps > 0$ and for $1 < p < + \infty$, they consider the functionals:
 $$\cF_\eps (u) := \int_\Omega f \left(\frac{x}{\eps}, u(x)\right)\; \de x, \qquad u \in L^p(\Omega; \R{d}), \; \cA u= 0,$$
\noindent where the integrand $f$ satisfies the following hypotheses:
\begin{itemize}
\item [(C1)] $f: \R{N}\times \R{d} \to [0, +\infty[$ is a continuous function $Q$-periodic in the first argument;
\item[(C2)] there exists $C > 0$ such that
$$ 0 \leq f(x, v) \leq C( 1 + |v|^p), \; \; \forall (x, v) \in \R{N}\times \R{d};$$
\item [(C3)] there exists $C >0$ such that
$$ f(x,v) \geq \frac{1}{C}|v|^p - C, \; \forall (x,v) \in \R{N}\times \R{d}.$$
\end{itemize}
Defining
$$ \Gamma-\liminf \cF_\eps (v) := \inf \left\{\liminf_{n\to +\infty} \cF_\eps(v_n) : v_n  \in L^p(\Omega; \R{d})\cap \; \ker\cA, \; v_n \wto v \; \text{in} \; L^p(\Omega; \R{d}) \right\},$$
and
$$ \Gamma-\limsup \cF_\eps (v) := \inf \left \{ \limsup_{n\to +\infty} \cF_\eps(v_n) : v_n  \in L^p(\Omega; \R{d})\cap \; \ker\cA, \; v_n \wto v \; \text{in} \; L^p(\Omega; \R{d}) \right\},$$

Their main result reads as follows:
\begin{theorem} \label{hombfl}
Under hypotheses (C1)--(C2), 
$$\mathcal{F}_{\hom} = \Gamma-\liminf \cF_\eps $$
where
$$\mathcal{F}_{\hom}(v; D) := \int_D f_{\hom}(v)\; \de x,$$
for all $v \in L^p(\Omega; \R{d}) \cap \ker\cA$ and $D \in \mathcal{O}(\Omega)$, with $f_{\hom}$ given by
$$ f_{\hom}(v) := \inf_{k \in \N{}} \frac{1}{k^{N}} \inf \left\{ \int_{kQ} f(x, v + w(x))\; \de x: w \in L^p_{kQ-\per}(\R{N}; \R{d})\cap\ker\cA, \; \int_{kQ} w(x)\; \de x = 0 \right\}.$$
Moreover, if also (C3) holds, 
$$\mathcal{F}_{\hom} = \Gamma-\lim  \cF_\eps.$$
\end{theorem}

\subsection{Case $p=1$}
Allowing for linear growth, i.e., taking $p=1$  implies working with sequences that are only bounded in $L^1$ and hence that can converge weakly-* (up to a subsequence) to some bounded Radon measure.

In \cite{ADMaso_92} , given a quasiconvex function $f: \R{N} \to \R{d}$ with linear growth, the authors derived an integral representation in  $BV(\Omega; \R{d})$ for the functional $\overline{F}$ arising from the relaxation in the $L^1_{\loc}(\Omega; \R{d})$ topology of the functional defined in $C^1(\Omega; \R{d})$ by
$$\int_\Omega f(\nabla u)\; \de x.$$
In \cite{FM1} the relaxation of functionals  of the form 
$$I(u) := \int_\Omega f(x, u(x), \nabla u(x))\; \de x, \qquad u \in BV(\Omega; \R{d}),$$
with $f(x, u,\cdot)$  quasiconvex with linear growth, was studied. 
This is a generalization of previous results, namely of the result obtained in \cite{DM} for the scalar case $d=1$.

In \cite{KriRin_CV_10} the relaxation of signed functionals with linear growth in the space BV was studied.

More precisely, for functionals defined for $u \in W^{1,1}(\Omega; \R{d})$ by $\int_\Omega f(\nabla u)\; \de x$, the authors derive strict continuity  and relaxation results in $BV(\Omega; \R{d}).$ The integrand $f: \R{d\times \N{}} \to \R{}$  is assumed to be  continuous, of linear growth at infinity and possibly unbounded from below. In \cite{rin1} a new proof of a  sequential weak-* lower semicontinuity result in $BV(\Omega; \R{d})$ of the form
$$ \mathcal{F}(u) := \int_\Omega f(x, \nabla u)\; \de x + \int_{S_u}f^\infty\left(x, \frac{D^s u}{|D^s u|}\right)\; \de|D^s u| + \int_{\partial \Omega} f^\infty (x, u \otimes \nu_\Omega)\; \de\mathcal{H}^{N-1}, \; u \in BV(\Omega; \R{d}),$$
was derived (here $\otimes$ denotes the tensor product). In contrast with the previous results, this proof uses a rigidity argument for gradients.
The generalization of this result in the context of $\mathcal{A}$-quasiconvexity was done in \cite{BCMS}, using rigidity arguments for functions in the kernel of the operator $\cA$.

In \cite{DEAGAR} a homogenization problem for functionals defined on vector-valued functions and with linear growth is addressed. 
The same problem in the scalar case was treated in \cite{Bou}, \cite{CEDA} and \cite{A}.
The extension to the vectorial case relies on the blow-up techniques introduced in \cite{ADMaso_92} and \cite{FM1}.
In \cite{MMS} a homogenization result in the context of $\cA$-quasiconvexity for integral functionals with linear growth was derived. This work extends to the case $p=1$ the homogenization results derived in \cite{BFL}.
The linear growth condition implies that concentration effects may appear and they need to be treated by carefully applying homogenization techniques in the setting of weak-* convergence in measure.

\smallskip
\par For $\Omega\subset\R{N}$, $N\geq2$,  a bounded open set, and, for $d\geq1$,
let $f:\Omega{\times}\R{d}\to[0,+\infty)$ be a non-negative measurable function in the first variable and Lipschitz continuous in the second, that satisfies the following linear growth-coercivity condition: there exist $C_1,C_2>0$ such that
\be\label{200}
C_1|\zeta|\leq f(x,\zeta)\leq C_2(1+|\zeta|)\qquad\text{for all $(x,\zeta)\in\Omega{\times}\R{d}$}.
\ee
Moreover $x\mapsto f(x,\zeta)$ is assumed to be $Q$-periodic for each $\zeta\in\R{d}$. 

\par The linear first order partial differential operator $\cA:\cD'(\Omega;\R{d})\to\cD'(\Omega;\R{M})$ 
is assumed to satisfy, in addition to  \eqref{murat}, the following condition
\be
\label {spanc}\mathrm{Span}(\cC)=\R{d},
\ee 
where $\cC$ stands for the characteristic cone associated with the operator $\cA$; see Definition \eqref{356}.

A representation theorem for the functional
\be \label{101}
\cF(\mu):=\inf\left\{\liminf_{n\to\infty} \int_\Omega f\left(\frac{x}{\eps_n},u_n(x)\right)\,\de x,\;\;\begin{array}{l}
\{u_n\}\subset\Lp1{}(\Omega;\R{d}),\;\;\cA u_n\stackrel{W^{-1,q}}{\longrightarrow}0, \\
u_n\wsto\mu,\,\, |u_n|\wsto\Lambda,\;\; \Lambda(\partial\Omega)=0
\end{array}\right\},
\ee
is derived, where $\cA$ is defined in \eqref{operator} and satisfies \eqref{murat}, \eqref{spanc}, and  $q\in(1,\frac{N}{N-1})$.
%



The main result of \cite{MMS} is contained in the following theorem.
\begin{theorem}\label{main}
Let $\Omega\subset\R{N}$ be a bounded open set, and 
let $f:\Omega{\times}\R{d}\to[0,+\infty)$ be a function which is measurable and $Q$-periodic in the first variable and Lipschitz continuous in the second, satisfying the growth condition \eqref{200}. 
For any $b\in\R{d}$, define 
\be\label{103}
f_{\cA-\hom}(b):=\inf_{R\in\N{}} \inf\left\{\ave_{RQ} f(x,b+w(x))\,\de x,\quad w\in\Lp1{RQ-\per}(\R{N};\R{d})\cap\ker\cA,\;\;\ave_{RQ} w=0\right\},
\ee
where $R\in\N{}$ and $f_{\cA-\hom}^\infty$ is the recession function of $f_{\cA-\hom}$.
For every $\mu\in\cM(\Omega;\R{d})\cap\ker\cA$, let $\mu=u^a\cL^N+\mu^s$ and let
$$
\cF_{\cA-\hom}(\mu):=\int_\Omega f_{\cA-\hom}(u^a)\,\de x+\int_\Omega f_{\cA-\hom}^\infty\left(\frac{\de\mu^s}{\de|\mu^s|}\right)\,\de|\mu^s|.
$$
Then, $\cF(\mu)=\cF_{\cA-\hom}(\mu)$.
\end{theorem}

\section{Remarks, conclusions, and perspective work}\label{sec:persp}

In this note, which we intend to be an introduction to the subject, we restricted ourselves to the case where just one microscopic scale is considered. 
This has the advantage of keeping the exposition (notationally) simple, yet deep enough for the reader to grasp the fundamental notions.
Since heterogeneities are small compared with global dimensions, usually different scales are used to describe the material: a macroscopic scale describes the behavior of the bulk, while at least one microscopic scale describes the heterogeneities of the composite. An important tool to address the case of two-scale homogenization was developed in the works of G.\@ Nguetseng \cite{NG} and G.\@ Allaire \cite{Allaire} (see also \cite {BA} and \cite{LNW}).

The theory was further developed by Allaire \cite{Allaire} by studying some general properties of two-scale convergence and applying it to several homogenization problems. 
Two-scale convergence has also been generalized to $n$-scale convergence in the obvious way (see \cite{BA}, \cite{LLPW1}, and \cite{LLPW2}). 
Moreover, it has been generalized to larger classes of admissible functions or measures. We refer to \cite{LNW} where the reader can find, in a self-contained way, the details and basic ideas of the theory.

Another important tool to address homogenization problems is the \emph{unfolding operator method} introduced by Cioranescu, Damlamian, and Griso in \cite{Dm2}. We refer the reader to \cite{D}, \cite{Dm1} and \cite{Dm3}.

\medskip

It is natural  to look for an extension of the  work in \cite{MMS} by weakening the conditions on $f$, namely the condition of Lipschitz continuity. 
The condition $\Lambda(\partial \Omega) = 0$ in \eqref{101} rules out concentrations on the boundary. 
Without this condition, the authors would have to restrict the class of operators $\cA$ by requiring appropriate extension properties. 
We refer the reader to \cite{KrRin} for some remarks about extension properties. 
A natural question that arises is precisely the characterization of operators $\cA$ that admit an extension to all of $\R{N}$.
Among other open problems, we would like to mention the case when $\cA$ depends on the space variable $x$, and the one where Murat's condition \eqref{murat} is not satisfied. 
We anticipate that a great effort must be spent for solving these problems, that are very challenging from the theoretical and technical points of view.
The solution of these problems will contribute to both completing the theoretical framework and providing powerful tools to study problems relevant for applications to engineering and physics.

\end{document}